\documentclass[a4paper,11pt]{article}
\usepackage{latexsym}
\usepackage{amsmath}
\usepackage{amssymb}
\usepackage{enumerate}
\usepackage{theorem}
\usepackage{array}
\usepackage[dvipdfmx]{graphicx}
\newtheorem{theorem}{Theorem}[section]
\newtheorem{lemma}[theorem]{Lemma}
\newtheorem{corollary}[theorem]{Corollary}

\theorembodyfont{\rmfamily}

\newtheorem{notation}[theorem]{Notation}

\newtheorem{remark}[theorem]{Remark}

\newcommand{\proof}{\noindent \mbox{\em Proof.\hspace*{2mm}}}
\newcommand{\qed}{\hfill \mbox{$  \Box $}}

\newcommand{\ssgyokan}{\vskip 8pt}

\title{The second largest number of points of plane curves over finite fields}
\author{
Masaaki Homma
\thanks{Partially supported by Grant-in-Aid
for Scientific Research (15K04829), JSPS.}
\\
 Department of Mathematics and Physics\\
Kanagawa University\\
Hiratsuka 259-1293, Japan\\
homma@kanagawa-u.ac.jp
\and
Seon Jeong Kim
\thanks{Partially supported by Basic Science Research Program through the National Research Foundation of Korea(NRF) 
funded by the Ministry of Education, Science and Technology (2012R1A1A2042228).
}\\
 Department of Mathematics and RINS\\
Gyeongsang National University\\
Jinju 660-701, Korea \\
skim@gnu.kr
}
\date{}

\begin{document}
\maketitle
\begin{abstract}
A basis of the ideal of the complement of a linear subspace
in a projective space over a finite field is given.
As an application, the second largest number of points of plane curves of degree $d$ over the finite field of $q$ elements is also given
for $d\geq q+1$.
\\
{\em Key Words}: Finite field, Basis of the ideal, Plane curve
\\
{\em MSC}: 11G20, 13F20, 14G15, 14N05
\end{abstract}

\section{Introduction}
Let $C$ be a curve of degree $d$ in projective plane $\mathbb{P}^2$ defined over the finite field  $\mathbb{F}_q$ of $q$ elements, which has no
$\mathbb{F}_q$-line components,
and $F(X, Y, Z)$ a homogeneous polynomial in $\mathbb{F}_q[X, Y, Z]$
which defines the curve $C$.
We are interesting in the set
\[
 C(\mathbb{F}_q) :=
    \{
      (a,b,c) \in \mathbb{P}^2(\mathbb{F}_q) \mid F(a,b,c)=0
    \},
\]
where $\mathbb{P}^2(\mathbb{F}_q)$ denotes the set of
$\mathbb{F}_q$-points of $\mathbb{P}^2$.
The number of $C(\mathbb{F}_q)$ is denoted by $N_q(C)$.

In the series of papers \cite{hom-kim2009, hom-kim2010a, hom-kim2010b},
we proved the Sziklai bound, which says that
\begin{quote}
$N_q(C) \leq (d-1)q +1$ unless $C$ is a curve over $\mathbb{F}_4$ which is
projectively equivalent to the curve defined by
\[
(X+Y+Z)^4 + (XY+YZ+ZX)^2+XYZ(X+Y+Z)=0
\]
over $\mathbb{F}_4$.
\end{quote}

In order to give a brief explanation of what we will do,
we should explain some notation.

\begin{notation}
The number of points of $\mathbb{P}^n(\mathbb{F}_q)$ is frequently denoted by
$\theta_q(n)$, that is
$\theta_q(n) = \frac{q^{n+1}-1}{q-1}$.

Let $x_0, \dots , x_n$ be coordinates of $\mathbb{P}^n$,
and $f_1, \dots , f_n$ homogeneous polynomials over $\mathbb{F}_q$.
The algebraic set in $\mathbb{P}^n$ over the algebraic closure of $\mathbb{F}_q$ defined by $f_1=\dots = f_n =0$ is frequently denoted by
$\{ f_1=\dots = f_n =0 \}.$

Here we summarize symbols related to plane curves, which will be used in this paper or future,
and also agree with ones in \cite{hom-kim2011}.

\begin{itemize}
    \item ${\cal C}_d(\mathbb{F}_q)$: the set of plane curves of degree $d$
over $\mathbb{F}_q$ without $\mathbb{F}_q$-linear components.
   \item ${\cal C}_d^i(\mathbb{F}_q) :=
       \{C \in {\cal C}_d(\mathbb{F}_q) \mid  
         \mbox{\rm $C$ is absolutely irreducible} \}$
   \item ${\cal C}_d^s(\mathbb{F}_q) :=
       \{C \in {\cal C}_d(\mathbb{F}_q) \mid  
         \mbox{\rm $C$ is nonsingular} \}$
   \item $M_q(d) = \max \{ N_q(C) \mid C \in {\cal C}_d(\mathbb{F}_q) \}$
   \item $M_q^i(d) = \max \{ N_q(C) \mid C \in {\cal C}_d^i(\mathbb{F}_q) \}$
   \item $M_q^s(d) = \max \{ N_q(C) \mid C \in {\cal C}_d^s(\mathbb{F}_q) \}$
   \item ${}_2M_q(d) = \max \{ N_q(C) \mid C \in {\cal C}_d(\mathbb{F}_q), \ 
               N_q(C)< M_q(d) \}$
   \item ${}_2M_q^i(d) = \max \{ N_q(C) \mid C \in {\cal C}_d^i(\mathbb{F}_q), \
               N_q(C)< M_q^i(d) \}$
   \item ${}_2M_q^s(d) = \max \{ N_q(C) \mid C \in {\cal C}_d^s(\mathbb{F}_q), \
               N_q(C)< M_q^s(d) \}$
 \end{itemize}
\end{notation}

Since the number of $\mathbb{P}^2(\mathbb{F}_q)$ is $\theta_q(2)$,
the Sziklai bound makes sense in the range $2 \leq d \leq q+2$.
We should add a few words to the Sziklai bound:
for $d =2$, $\sqrt{q}+1$ (if $q$ is square), $q-1$, $q$, $q+1$ and $q+2$,
this bound is sharp.
For those degrees $d$,
it might be interesting to know the second largest number
${}_2M_q(d)$.
In this paper, we try to find it for $d= q+1$ and $q+2$.
Furthermore, since we have already knew $M_q(d)= \theta_q(2)$ if
$d \geq q+3$ \cite[Prop. 1.1]{hom-kim2009},
we can determine ${}_2M_q(d)$ for $d \geq q+3$.

More precisely, we will show the following facts.

\begin{theorem}\label{maintheorem}
 \begin{enumerate}[\rm (1)]
    \item Suppose that $q>3$.
  \begin{enumerate}[\rm (i)]
  \item If $d \geq 2q-1$, then 
     ${}_2M_q(d) = \theta_q(2) -1 = q^2+q.$
  \item If $2q-1 \geq d \geq q+2$, then
     ${}_2M_q(d) = \theta_q(2) -(2q-d) = q^2 + d - q +1.$
  \item If $d=q+2$, then
    ${}_2M_q(q+2) ={}_2M_q^i(q+2) =q^2 + 3.$
  \end{enumerate}
   \item For any $q$, we have
$
 {}_2M_q(q+1) ={}_2M_q^i(q+1) =q^2.
$
 \end{enumerate}
\end{theorem}

\begin{figure}
\centering
\includegraphics[width=16cm]{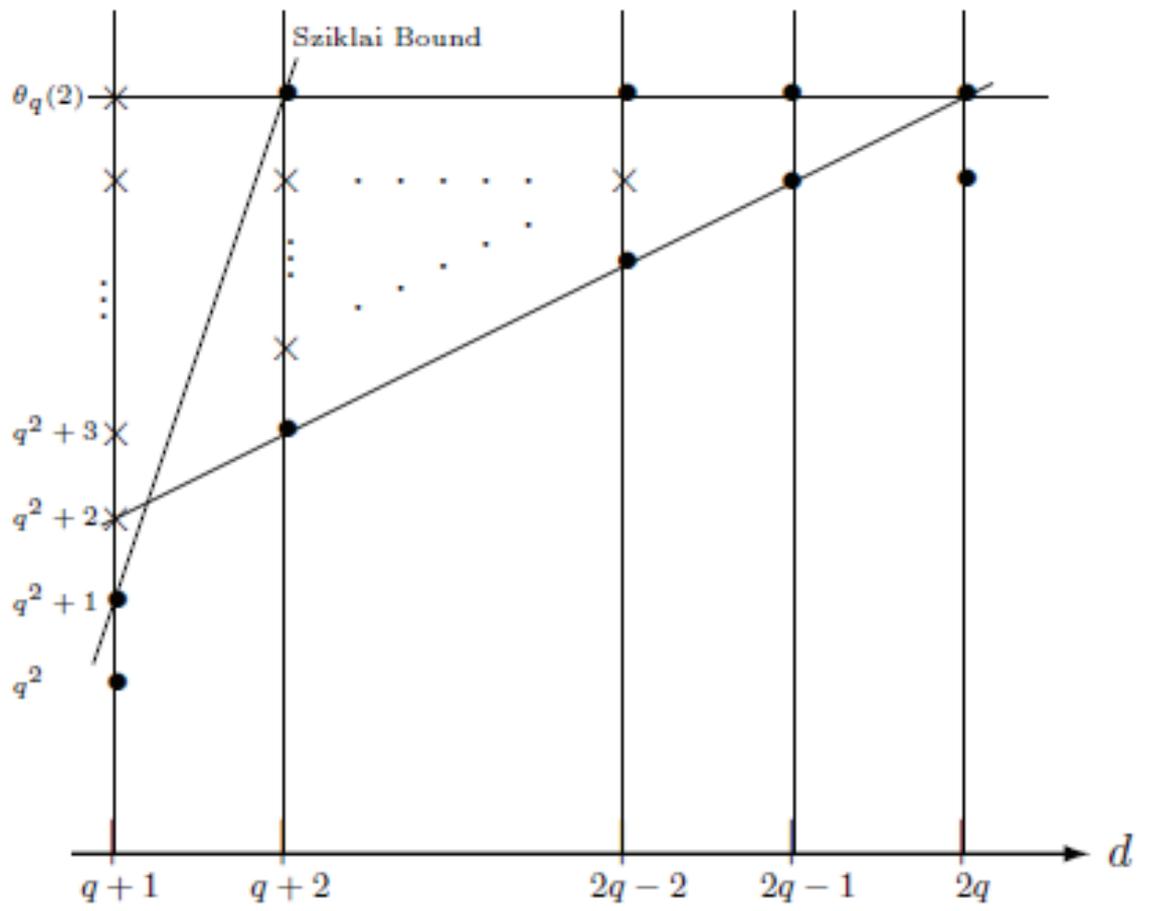}\\
\caption{$M_q(d)$ and ${}_2M_q(d)$ for $d \geq q+1$ with $q>3$}
 \label{picture1}
\end{figure}

Figure~\ref{picture1} illustrates the main theorem.
In the picture $\bullet$ indicates existence of a plane curve
of assigned degree and number of points,
and $\times$ non-existence.

\vskip 4pt

The important ingredient of the proof is to determine
a basis of the ideal of
$\mathbb{P}^2(\mathbb{F}_q) \setminus \{ \text{an $\mathbb{F}_q$-point} \}$
for (1) of Theorem~\ref{maintheorem},
and $\mathbb{P}^2(\mathbb{F}_q) \setminus \{ \text{an $\mathbb{F}_q$-line} \}$
for (2).
We formulate them under a little more general setting
in Theorem~\ref{ideal} below.

\section{Ideal of the complement of a linear subspace in $\mathbb{P}^n(\mathbb{F}_q)$}
In this section, we consider two subsets of $\mathbb{F}_q$-points of
$\mathbb{P}^n$ with homogeneous coordinates $x_0, x_1, \dots , x_n$.

The following lemma, which is easy to show, is a kind of folklore.
\begin{lemma}
The ideal of $\mathbb{P}^n(\mathbb{F}_q)$
in $\mathbb{F}_q[x_0, \dots , x_n]$ is generated by
$
\{
x_i^q x_j - x_i x_j^q \mid 0 \leq i < j \leq n
\}.
$
\end{lemma}

\begin{theorem}\label{ideal}
Let $k$ be an integer with $1 \leq k \leq n$,
and $\mathbb{P}^{k-1}$ the linear subspace
defined by
$ x_k = x_{k+1} = \dots = x_n=0.$
Then the ideal of 
$\mathbb{P}^n(\mathbb{F}_q) \setminus \mathbb{P}^{k-1}(\mathbb{F}_q) $
in $\mathbb{F}_q[x_0, \dots , x_n]$ is generated by
\[
\{
x_i^q x_j - x_i x_j^q \mid 0 \leq i < j \leq n
\} \cup
\{
x_s \prod_{i=k}^n (x_i^{q-1} - x_s^{q-1})
\mid s = 0, 1, \dots , k-1
\}.
\]
\end{theorem}
\proof
(Step 1)
In the first step, we handle the case
where the linear subspace is of dimension $0$,
that is, $k=1$.

Let $J_n$ be the ideal of $\mathbb{P}^n(\mathbb{F}_q) \setminus \{ (1,0, \dots , 0) \}$ in $\mathbb{F}_q[x_0, \dots , x_n]$,
and $I_n$ the ideal generated by
\[
\{
x_i^q x_j - x_i x_j^q \mid 0 \leq i < j \leq n
\} \cup
\{
x_0 \prod_{i=1}^n (x_i^{q-1} - x_0^{q-1})
\}.
\]
It is obvious that $J_n \supset I_n$.
We will show $J_n = I_n$ by induction on $n$.
For $n=1$,
\[
I_1 =(x_0^q x_1 - x_0 x_1^q, x_0(x_1^{q-1}- x_0^{q-1}) )
     = (x_0(x_1^{q-1}- x_0^{q-1}) )
\]
is trivially the ideal of $\{ (\lambda , 1) \mid \lambda \in \mathbb{F}_q \}
= \mathbb{P}^1(\mathbb{F}_q) \setminus \{ (1,0) \}$.

Next we assume that $n \geq 2$.
Let $F(x_0, \dots , x_n)$ be a homogeneous polynomial of degree $d$ in $J_n$.
We want to show that $F(x_0, \dots , x_n) \equiv 0 \mod I_n$.
First we prove the following claim:

{\it Claim.\ } We have
\[
F(x_0, \dots , x_n) \equiv \sum_{i=0}^{q-1} f_i(x_0, \dots , x_{n-1})x_n^i 
\mod I_n,
\]
where $f_i(x_0, \dots , x_{n-1})$ is homogeneous of degree $d-i$.

When $d \leq q-1$, there is nothing to do.
Suppose $d \geq q$.
For a monomial $M=x_0^{e_0}\cdots x_n^{e_n}$ which appears in $F$,
if $q \leq e_n<d$, then
$
M \equiv  x_0^{e'_0}\cdots x_{n-1}^{e'_{n-1}} x_n^{e'_n} \mod I_n
$
with $e'_n < e_n$,
because $x_i^q x_n - x_i x_n^q \in I_n$.
(Note that since $e_n <d$, there is an $e_i$ with $0 \leq i \leq n-1$ such that $e_i \geq 1$.)
Repeating this procedure, we get
$M \equiv x_0^{e'_0}\cdots x_{n-1}^{e'_{n-1}} x_n^{e'_n} \mod I_n$
with  $e'_n < q$.
Hence
\begin{equation}\label{fmodin}
F(x_0, \dots , x_n) \equiv cx_n^d +
\sum_{i=0}^{q-1} f_i(x_0, \dots , x_{n-1})x_n^i 
\mod I_n.
\end{equation}
Since the right-hand side of (\ref{fmodin}) is also an element of $J_n$,
$(0, \dots, 0, 1)$ is its zero. Hence $c=0$.

From this claim,
we may assume that
\[
F(x_0, \dots , x_n) = \sum_{i=0}^{q-1} f_i(x_0, \dots , x_{n-1})x_n^i , \
 \deg f_i = d-i.
\]
Let $(\lambda_0, \dots , \lambda_{n-1}) \in
\mathbb{F}_q^n \setminus \{(1,0, \dots , 0) \}$.
Since
$F(\lambda_0, \dots , \lambda_{n-1}, \mu)=0$
for any $\mu \in \mathbb{F}_q$
and $F(\lambda_0, \dots , \lambda_{n-1}, x_n)$ is of degree $q-1$ in $x_n$,
$f_{\nu}(\lambda_0, \dots , \lambda_{n-1})=0$ for $\nu =0, \dots , q-1$.
Therefore $f_{\nu}(x_0, \dots , x_{n-1}) \in I_{n-1}$
by the induction hypothesis,
that is, 
\[
f_{\nu}(x_0, \dots , x_{n-1}) =
  \sum_{0 \leq i <j\leq n-1} h_{ij}\cdot(x_i^qx_j - x_ix_j^q)
   + k_{\nu} x_0\prod_{i=1}^{n-1} (x_i^{q-1}-x_0^{q-1}),
\]
where $k_{\nu} = k_{\nu}(x_0, \dots , x_{n-1})$
is of degree $d - \nu -((n-1)(q-1)+1)$, say $d_{\nu}$.
(Note that if $d - \nu -((n-1)(q-1)+1) <0$, we understand $k_{\nu}=0$.)
Put
\[
k_{\nu}(x_0, \dots , x_{n-1}) = c_{\nu}x_0^{d_\nu}
   + k_{\nu}^l(x_0, \dots , x_{n-1}),
\]
where the degree in $x_0$ of $k_{\nu}^l$ is less than $d_{\nu}$.
Then
\[
F(x_0, \dots , x_n) \equiv 
\sum_{\nu=0}^{q-1} (c_{\nu}x_0^{d_\nu} + k_{\nu}^l(x_0, \dots , x_{n-1}))
x_0\prod_{i=1}^{n-1} (x_i^{q-1}-x_0^{q-1})x_n^{\nu} \ 
\mod I_n.
\]
So we may assume that
\[
F(x_0, \dots , x_n) = G(x_0, \dots , x_n) + H(x_0, \dots , x_n)
\]
with
\[
G(x_0, \dots , x_n) = x_0\prod_{i=1}^{n-1} (x_i^{q-1}-x_0^{q-1})
                             (\sum_{\nu=0}^{q-1}c_{\nu}x_0^{d_\nu}x_n^{\nu})
\]
and
\[
H(x_0, \dots , x_n) = \sum_{\nu=0}^{q-1}k_{\nu}^l(x_0, \dots , x_{n-1})
x_0\prod_{i=1}^{n-1} (x_i^{q-1}-x_0^{q-1})x_n^{\nu}.
\]
Since $k_{\nu}^l(x_0, \dots , x_{n-1})$ is the sum of lower terms in $x_0$
of the homogeneous polynomial $k_{\nu}(x_0, \dots , x_{n-1})$,
each term of $k_{\nu}^l$ contains some $x_j$ with $ 0< j \leq n-1$.
Since
\[
x_j x_0 \prod_{i=1}^{n-1} (x_i^{q-1}-x_0^{q-1})
    =\left( \prod^{n-1}_{
                  \begin{subarray}{c}
                  i=1\\
                  i \neq j
                  \end{subarray}
                   } 
                   (x_i^{q-1}-x_0^{q-1}) \right)
                   (x_j^{q}x_0 - x_j x_0^q) \in I_n,
\]
we have $H(x_0, \dots , x_n) \in I_n$.

Finally we show that $G(x_0, \dots , x_n) \in I_n$.
For $(\mu, \rho) \in \mathbb{F}_q^2 \setminus \{(1,0)\}$,
$G(\mu, 0, \dots, 0, \rho)=0$
because $G= F -H \in J_n$.
Hence $G(x_0, 0 , \dots , 0, x_n) \in I_1$ with respect to coordinates
$x_0, x_n$
by the induction hypothesis,
that is $G(x_0, 0 , \dots , 0, x_n) \in 
       \left( x_0(x_n^{q-1}-x_0^{q-1}) \right)$.
On the other hand, since
\[
G(x_0, 0 , \dots , 0, x_n) =(-1)^{n-1}x_0^{(n-1)(q-1)+1}
     (\sum_{\nu=0}^{q-1}c_{\nu}x_0^{d_\nu}x_n^{\nu}),
\]
there is a polynomial $u(x_0, x_n)$
such that
\[
\sum_{\nu=0}^{q-1}c_{\nu}x_0^{d_\nu}x_n^{\nu}=
    u(x_0, x_n)(x_n^{q-1}-x_0^{q-1}).
\]
Therefore
$G(x_0, \dots , x_n) = u(x_0, x_n)x_0\prod_{i=1}^{n} (x_i^{q-1}-x_0^{q-1})
\in I_n.$

(Step 2)
We complete the proof by induction on the dimension of the linear subspace
$\mathbb{P}^{k-1}$.
Suppose $k \geq 2$.
Let $J_{n,k}$ be the ideal of
$\mathbb{P}^n(\mathbb{F}_q) \setminus \mathbb{P}^{k-1}(\mathbb{F}_q) $,
and $I_{n,k}$ the ideal generated by
$
\{
x_i^q x_j - x_i x_j^q \mid 0 \leq i < j \leq n
\} \cup
\{
x_s \prod_{i=k}^n (x_i^{q-1} - x_s^{q-1})
\mid s = 0, 1, \dots , k-1
\}.
$
Note that the second subscript in notation of those ideals
indicates the dimension of the vector space
$H^0(\mathbb{P}^{k-1}, \mathcal{O}(1))$.
Obviously $J_{n,k} \supset I_{n, k}$.
Let $F(x_0, \dots , x_n)$ be a homogeneous polynomial in $J_{n,k}$.
Divide the terms of $F$ into two parts:
let $G_1(x_0, x_2, \dots , x_n)$ be the sum of terms of $F$ each of which
does not contains $x_1$, and $G_2(x_0, x_1, \dots , x_n)$
the sum of terms of $F$ each of which
contains $x_1$.
In other words,
\begin{align}
G_1(x_0, x_2, \dots , x_n) &= F(x_0, 0, x_2, \dots , x_n)  
\label{conditiononG1} \\
G_2(x_0, x_1, \dots , x_n) &= F(x_0, x_1, x_2, \dots , x_n) \notag
                        - F(x_0, 0, x_2, \dots , x_n).
\end{align}
From (\ref{conditiononG1}),
$G_1(\lambda_0, \lambda_2, \dots , \lambda_n) =0$
for any
$(\lambda_0, \lambda_2, \dots , \lambda_n) \in 
 \mathbb{P}^{n-1} \setminus \{x_k = \dots =x_n =0 \}$
with respect to the coordinates
$x_0, \Hat{x_1}, x_2, \dots, x_{k-1}; x_k,  \dots , x_n$, 
where $\Hat{x_1}$ means omitting $x_1$.
Hence, by induction hypothesis\footnote{
The dimension of the linear subspace $\{x_k = \dots =x_n =0 \}$
in $\mathbb{P}^{n-1}$ with respect to the coordinates
$x_0, \Hat{x_1}, x_2, \dots , x_n$ is $k-2$.
},
$G_1(x_0, x_2, \dots , x_n) \in I_{n-1, k-1}$
with respect to the coordinates
$x_0, \Hat{x_1}, x_2, \dots, x_{k-1}; x_k,\dots , x_n$,
that is,
$G_1$ is contained in the ideal generated by
\[
\{
x_i^q x_j - x_i x_j^q \mid 0 \leq i < j \leq n, \ 
   i, j  \neq 1
\} \cup
\{
x_s \prod_{i=k}^n (x_i^{q-1} - x_s^{q-1})
\mid s = 0, 2, \dots , k-1
\}.
\]
Hence $G_1$ is an element of the original
$I_{n, k}$,
and hence $G_2 = F - G_1 \in J_{n,k}$.
Since $x_0^qx_1 - x_0x_1^q \in I_{n,k}$
and each term of $G_2(x_0, x_1, \dots , x_n)$
contains $x_1$,
\begin{equation}\label{reduce}
G_2(x_0, x_1, \dots , x_n) \equiv 
   \sum_{i=0}^{q-1} g_i(x_1, \dots , x_n) x_0^i
      \mod I_{n,k}.
\end{equation}
Hence it is enough to observe that the right-hand side of
(\ref{reduce}) belongs $I_{n,k}$.
Since the polynomial already contained in $J_{n,k}$,
for a fixed $(\lambda_1, \dots , \lambda_n) \in
\mathbb{F}_q^n \setminus \{ \lambda_k = \dots = \lambda_n=0 \}$
and an arbitrary $\mu \in \mathbb{F}_q$,
\[
\sum_{i=0}^{q-1} g_i(\lambda_1, \dots , \lambda_n)\mu^i =0.
\]
Since $\sum_{i=0}^{q-1} g_i(\lambda_1, \dots , \lambda_n) x_0^i$ 
is a polynomial of degree $q-1$ in $x_0$,
we have $ g_i(\lambda_1, \dots , \lambda_n)=0$
for each $i = 0, 1, \dots , q-1$.
Hence $g_i(x_1, \dots , x_n) \in I_{n-1, k-1}$
with respect to the variables
$x_1, \dots x_{k-1};x_k, \dots , x_n$
by induction hypothesis.
It is obvious that this ideal $I_{n-1, k-1}$ is contained in $I_{n, k}$.
\qed

\ssgyokan
The following corollary is just the case $k=1$ in Theorem~\ref{ideal},
which has been proved in Step 1.
\begin{corollary}\label{ideallessonepoint}
The ideal of $\mathbb{P}^n(\mathbb{F}_q) \setminus \{ (1,0, \dots , 0) \}$
in $\mathbb{F}_q[x_0, \dots , x_n]$ is generated by
\[
\{
x_i^q x_j - x_i x_j^q \mid 0 \leq i < j \leq n
\} \cup
\{
x_0 \prod_{i=1}^n (x_i^{q-1} - x_0^{q-1})
\}.
\]
\end{corollary}

\begin{corollary}\label{minimumdegree}
Let $P_0 \in \mathbb{P}^n(\mathbb{F}_q)$.
Then there is
a homogeneous polynomial $F(x_0,\dots , x_n)$ of degree $d$
in $\mathbb{F}_q[x_0, \dots , x_n]$
such that 
the hypersurface $H$  defined by $F=0$
satisfies $H(\mathbb{F}_q) = \mathbb{P}^n(\mathbb{F}_q) \setminus \{ P_0 \}$
if and only if $d \geq (q-1)n +1$.
\end{corollary}
\proof
It is obvious that we may choose $P_0= (1, 0, \dots, 0)$.
Since polynomials $x_i^q x_j - x_i x_j^q$ ($0 \leq i < j \leq n$)
vanish on the entire $\mathbb{P}^n(\mathbb{F}_q)$,
the ``only if" part follows from the above corollary immediately.
Conversely the hypersurface $H$ defined by
$x_0^{d-(q-1)n} \prod_{i=1}^n (x_i^{q-1} - x_0^{q-1})$
has the desired property if $d>(q-1)n$.
\qed

\begin{corollary}\label{affineideal}
 Let $\mathbb{A}^n$ be the affine part of $\mathbb{P}^n$
defined by $x_n \neq 0$.
The ideal of $\mathbb{A}^n(\mathbb{F}_q)$ in $\mathbb{F}_q[x_0, \dots , x_n]$
is generated by
\[
\{ x_sx_n^{q-1} -  x_s^q \mid 0 \leq s < n \}.
\]
\end{corollary}
\proof
From Theorem~\ref{ideal},
the ideal is generated by
\[
\{
x_i^q x_j - x_i x_j^q \mid 0 \leq i < j \leq n
\}
\cup
\{ x_sx_n^{q-1} -  x_s^q \mid 0 \leq s < n \}.
\]
Since
\[
 x_i^qx_j - x_ix_j^q 
  =(x_jx_n^{q-1} - x_j^q)x_i
        - (x_ix_n^{q-1} - x_i^q)x_j,
\]
the ideal is generated by only the latter set.
\qed

\section{The second largest number for $d \geq q+1$}
Throughout this section we fix projective plane
$\mathbb{P}^2$ with coordinates $X, Y, Z$ over $\mathbb{F}_q$.

The following lemma is a corollary of the Sziklai bound,
which guarantees a plane curve without $\mathbb{F}_q$-line components
having many $\mathbb{F}_q$-points to be absolutely irreducible.
\begin{lemma}\label{absolutelyirreducible}
Let $C$ be a curve over $\mathbb{F}_q$ of degree $d$ in $\mathbb{P}^2$
without $\mathbb{F}_q$-linear components.
If $N_q(C) \geq (d-2)q +3$, then $C$ is absolutely irreducible.
\end{lemma}
\proof
See \cite[Corollary~2.2]{hom-kim2015preprint}.
\qed

\subsection{The case $q+2 \leq d \leq 2q-1$}
\begin{lemma}\label{mainlemma}
Let $C$ be a curve over $\mathbb{F}_q$ of degree $d$ in $\mathbb{P}^2$
with $q+2 \leq d \leq 2q-2$,
which may have $\mathbb{F}_q$-linear components.
Suppose $C(\mathbb{F}_q) \neq \mathbb{P}^2(\mathbb{F}_q)$.
Let
$N_q(C)= \theta_q(2) - (m+1)$.
Then $m+1 \geq 2q-d$, that is $N_q(C) \leq q^2 +(d-q+1)$;
and if equality occurs, then those $m+1$ missing points are collinear. 
\end{lemma}
\proof
Let $F \in \mathbb{F}_q[X, Y, Z]$ be an equation of $C$.
Put
$\mathbb{P}^2(\mathbb{F}_q) \setminus C(\mathbb{F}_q)
  = \{P_0, \dots , P_m \}$.
For each $i$ with $1 \leq i \leq m$,
choose a linear form $L_i$ over $\mathbb{F}_q$
such that the line $\{ L_i = 0 \}$
passes through $P_i$ but does not $P_0$.
Then
$\{ L_1 \cdots L_mF = 0\}(\mathbb{F}_q) = 
\mathbb{P}^2(\mathbb{F}_q) \setminus \{P_0\}$.
Hence $\deg  L_1 \cdots L_mF \geq 2q-1$ by Corollary~\ref{minimumdegree}.
Therefore $m+1 \geq 2q-d$.

Next we will show that $P_0, \dots , P_m $ are collinear if $m=2q-d-1$.
For two indexes $i$, $j$ with $1 \leq i < j \leq m$,
if the line $\overline{P_iP_j}$ joining $P_i$ and $P_j$
does not contain $P_0$, then
the curve
\[
D:= \overline{P_iP_j} \cup \{ \prod_{
                            \begin{subarray}{c}
                                k  \text{ with }\\
                                 k \neq i, j
                            \end{subarray}
                                 }L_k \cdot F =0 \}
\]
is of degree 
$
m-1 +d= 2q-2
$
and $D(\mathbb{F}_q) = \mathbb{P}^2(\mathbb{F}_q) \setminus \{P_0\}$,
which contradicts to Corollary~\ref{minimumdegree}.
Hence $P_0, \dots , P_{m} $ are collinear.
\qed
\begin{theorem}\label{theoremforbigd}
Fix an integer $d$ with $q+2 \leq d \leq 2q-1$.
Let $C$ be a plane curve over $\mathbb{F}_q$ of degree $d$
such that $N_q(C)=q^2+(d-q+1)$.
Then $C$ is projectively equivalent to the curve defined by
\[
(X^q-XZ^{q-1})f(X,Y,Z) + (Y^q-YZ^{q-1})g(X,Y,Z)=0
\]
over $\mathbb{F}_q$, where 
$f(X,Y,Z)$ and $g(X,Y,Z)$ are of degree $d-q$ and
$Xf(X,Y,0) + Yg(X,Y,0)= 0$
has $d-q+1$ distinct roots on $\mathbb{P}^1(\mathbb{F}_q)$
with homogeneous coordinates $X, Y$.

Conversely such a curve has $q^2+(d-q+1)$ $\mathbb{F}_q$-points.
{\rm(}But we don't know whether the curve has 
an $\mathbb{F}_q$-line component or not.{\rm)}
Moreover if $q>3$, then
there is a curve without $\mathbb{F}_q$-line components
among those curves.
\end{theorem}
\proof
From Lemma~\ref{mainlemma},
the missing $2q-d$ points are collinear.
Choosing coordinates $X, Y, Z$ of $\mathbb{P}^2$
as those points lie on the line $Z=0$.
Hence the equation of $C$ is in the ideal
$
(X^q - XZ^{q-1}, Y^q - YZ^{q-1})
$
by Corollary~\ref{affineideal}.
So $C$ is defined by
\[
(X^q - XZ^{q-1})f(X,Y,Z) + (Y^q - YZ^{q-1})g(X,Y,Z)=0,
\]
where $f(X,Y,Z), g(X,Y,Z) \in \mathbb{F}_q[X, Y, Z]$
are homogeneous of degree $d-q$.
Since $C(\mathbb{F}_q) \cap \{ Z=0 \}$ consists of
$d-q+1$ points,
\begin{equation}\label{condition1}
 X^q f(X,Y,0) + Y^qg(X,Y,0)=0
\end{equation}
has $d-q+1$ roots on $\mathbb{P}^2(\mathbb{F}_q)$
with coordinates $X, Y$.
Since we consider only $\mathbb{F}_q$-points,
(\ref{condition1}) is equivalent to
\begin{equation}\label{condition2}
X f(X,Y,0) + Y g(X,Y,0)=0.
\end{equation}

Conversely 
it is obvious that such a curve has
$q^2+(d-q+1)$ $\mathbb{F}_q$-points.

Now we construct a required equation if $q>3$.
Since $d \leq 2q-1$, we can choose $d-q+1$ distinct
elements $\alpha_1, \dots , \alpha_{d-q+1}$ in $\mathbb{F}_q$.
Let $\beta_0 (=1),  \beta_1 \dots , \beta_{d-q+1}$
be elements of $\mathbb{F}_q$
determined by
\[
\prod_{i=1}^{d-q+1} (Y- \alpha_i X)
   = \sum_{i=0}^{d-q+1} \beta_{d-q+1-i} X^{d-q+1-i} Y^i.
\]
For $\mathbf{c} = (c_1, \dots , c_{d-q}) \in \mathbb{F}_q^{d-q}$,
we consider the polynomial
\begin{multline}\label{fc}
F_{\mathbf{c}}(X,Y,Z) = \\
(X^q - XZ^{q-1})( \sum_{i=0}^{d-q} \beta_{d-q+1-i}X^{d-q-i} Y^i 
 + \sum_{i=1}^{d-q} c_i X^{d-q-i}Z^i ) \\
     +(Y^q - YZ^{q-1})Y^{d-q}.
\end{multline}
Let $C$ be the curve $F_{\mathbf{c}}(X,Y,Z) =0$.
The claims are
\begin{enumerate}[{Claim} 1.\ ]
 \item $N_q(C) = q^2 + (d-q+1)$;
 \item $C$ has no $\mathbb{F}_q$-line components if one choose
 $\mathbf{c}\in \mathbb{F}_q^{d-q} $ appropriately.
\end{enumerate}
For Claim 1, it is obvious $C(\mathbb{F}_q) \supset \mathbb{A}^2(\mathbb{F}_q)$, where $\mathbb{A}^2 = \mathbb{P}^2 \setminus \{ Z=0 \}$.
From (\ref{condition2}),
$C(\mathbb{F}_q) \cap \{ Z=0\}$ comes from
\begin{align}
 0 &= \sum_{i=0}^{d-q} \beta_{d-q+1-i} X^{d-q+1-i}Y^i + Y^{d-q+1} \notag \\
   &= \prod_{i=1}^{d-q+1} (Y- \alpha_i X). \notag
\end{align}
Hence $C(\mathbb{F}_q) \cap \{ Z=0\}
=\{ (1, \alpha_i , 0) \mid i= 1, 2, \dots , d-q+1\}$,
and $N_q(C) = q^2 + (d-q+1)$.

We see the second claim.
An $\mathbb{F}_q$-line which is a component of $C$ must pass through
one of the $d-q+1$ $\mathbb{F}_q$-points on $Z=0$ above,
and can't be the line $ Z=0$.
So the possible $\mathbb{F}_q$-lines that may be components of $C$
are
\[
\{Y - \alpha X = \rho Z \mid 
   \alpha \in \{\alpha_1, \dots , \alpha_{d-q+1} \} , \ 
                            \rho \in \mathbb{F}_q  \}.
\]
By direct computation,
$F_{\mathbf{c}}(X,\alpha X +\rho Z,Z)$ can be written as
\[
(X^q - XZ^{q-1})
 ( \sum_{i=0}^{d-q} \gamma_{i}(\alpha, \rho) X^{d-q-i}Z^i
       + \sum_{i=1}^{d-q} c_i X^{d-q-i}Z^i ).
\]
The number of vectors in
\[
\{
(\gamma_{0}(\alpha, \rho), \dots , \gamma_{d-q}(\alpha, \rho) ) \mid
 \alpha \in \{\alpha_1, \dots , \alpha_{d-q+1} \} , \ 
                            \rho \in \mathbb{F}_q  
\}
\]
is at most $(d-q+1)q$
and that of the vectors $\{(0, \mathbf{c})
         \mid \mathbf{c} \in  \mathbb{F}_q^{d-q} \}$ is $q^{d-q}$.
Since $(d-q+1)q < q^{d-q}$ if $q+2 \leq d \leq 2q-1$ and $q>3$,
we can choose a vector $\mathbf{c}$ such that
$F_{\mathbf{c}}(X,\alpha X + \rho Z,Z)$ is nontrivial
for all pairs $(\alpha , \rho)$.
\qed

\begin{remark}\label{justimitation}
Suppose that $q>3$ and $d \geq 2q$.
Let $\{ e_{\alpha} \mid \alpha \in \mathbb{F}_q \}$
be integers such that $ e_{\alpha} \geq 1$ for any $\alpha \in \mathbb{F}_q$
and $\sum_{\alpha \in \mathbb{F}_q} e_{\alpha} = d-q+1.$
Let $\beta_{0}, \dots , \beta_{d-q+1}$ be elements of $\mathbb{F}_q$
determined by 
\[
\prod_{\alpha \in \mathbb{F}_q} (Y- \alpha_i X)^{e_{\alpha}}
   = \sum_{i=0}^{d-q+1} \beta_{d-q+1-i} X^{d-q+1-i} Y^i
\]
as in the proof of Theorem~\ref{theoremforbigd}.
Then, by just imitating the proof of the claims,
we can find $(c_1, \dots , c_{d-g}) \in \mathbb{F}_q^{d-q}$
such that the curve $C$ with the equation $F_{\mathbf{c}}(X,Y,Z) =0$
defined by (\ref{fc})
has the property
\begin{enumerate}[(i)]
 \item $N_q(C) = \theta_q(2) -1$;
 \item no $\mathbb{F}_q$-line is a component of $C$.
\end{enumerate}
\end{remark}

\vskip 4pt

\noindent
\mbox{\em Proof of Theorem~{\rm \ref{maintheorem} (1).}\hspace*{2mm}}
(i) and (ii) are consequences of Theorem~\ref{theoremforbigd}
and Remark~\ref{justimitation}.
If $d=q+2$, any curve $C$ over $\mathbb{F}_q$
without $\mathbb{F}_q$-line components having $q^2+3$ $\mathbb{F}_q$-points
is absolutely irreducible by Lemma~\ref{absolutelyirreducible}.
\qed

\subsection{The case $d=q+1$}

\begin{theorem}
Let $\mathbb{A}^2$ be the affine part of $\mathbb{P}^2$ defined by $Z\neq 0$,
and $C$ a curve of degree $d$ in $\mathbb{P}^2$ defined over $\mathbb{F}_q$.

\begin{enumerate}[\rm (i)]
\item If $C(\mathbb{F}_q)=\mathbb{A}^2(\mathbb{F}_q)$, then $d \geq q+1$.
\item When $d=q+1$, $C(\mathbb{F}_q)=\mathbb{A}^2(\mathbb{F}_q)$
if and only if $C$ is defined by an equation of the following type:
\begin{equation}\label{equationqplusone}
(X^q-XZ^{q-1}, Y^q-YZ^{q-1})
\left(
 \begin{array}{ccc}
   a_0 & a_1 & a_2 \\
   b_0 & b_1 & b_2
 \end{array}
\right)
\left(
\begin{array}{c}
X \\ Y \\ Z
\end{array}
\right) = 0,
\end{equation}
where $a_0, \dots , b_2 \in \mathbb{F}_q$
with the polynomial
\begin{equation}\label{additionalcondition}
(s,t)
\left(
 \begin{array}{cc}
   a_0 & a_1 \\
   b_0 & b_1 
 \end{array}
\right)
\left(
\begin{array}{c}
s \\ t
\end{array}
\right)
\end{equation}
in $s$ and $t$ being irreducible over $\mathbb{F}_q$.
\item The curve described above is absolutely irreducible with
a unique singular points, which is an $\mathbb{F}_q$-points.
\end{enumerate}
\end{theorem}
\proof
(i)
Since $C(\mathbb{F}_q) = \mathbb{A}^2(\mathbb{F}_q)$,
the curve has no $\mathbb{F}_q$-linear component.
Hence we can apply the Sziklai bound to $C$,
namely, $q^2 \leq (d-1)q +1$.
So $d \geq q+1$.

(ii) 
From Corollary~\ref{affineideal},
the equation of 
$C$ can be written as (\ref{equationqplusone}).
Since the line $Z=0$ does not meet with $C(\mathbb{F}_q)$,
the polynomial (\ref{additionalcondition}) must be irreducible
over $\mathbb{F}_q$,
and vice versa.

(iii)
Since $\deg C= q+1$ and $N_q(C) = q^2$,
$C$ is absolutely irreducible if $q \geq 3$
by Lemma~\ref{absolutelyirreducible}.
Suppose $q=2$. Then $\deg C =3$.
Further suppose that $C$ is reducible.
Since $C$ can't have $\mathbb{F}_q$-line components,
$C$ must be a union of three conjugate  $\mathbb{F}_{q^3}$-lines
over $\mathbb{F}_q$.
Hence $N_q(C) \leq 1$, which is a contradiction.
Hence $C$ is absolutely irreducible even if $q=2$.

Since (\ref{additionalcondition}) is irreducible over $\mathbb{F}_q$,
$
 \det \left(
        \begin{array}{cc}
         a_0 & a_1 \\
         b_0 & b_1 
        \end{array}
      \right) \neq 0.
$
Let $(x_0, y_0)$ be the unique solution of 
\[
   \left(
        \begin{array}{cc}
         a_0 & a_1 \\
         b_0 & b_1 
        \end{array}
    \right)
    \left(
       \begin{array}{c}
       x \\ y
       \end{array}
    \right)
        = -
    \left(
       \begin{array}{c}
       a_2 \\ b_2
       \end{array}
    \right).
\]
Since
$(x_0, y_0, 1) \in \mathbb{A}^2(\mathbb{F}_q)$,
it is a point of $C$.
Choose new coordinates $X', Y', Z'$ as
\[
\left\{
 \begin{array}{ccl}
  X' & = & X-x_0Z \\
  Y' & = & Y-y_0Z \\
  Z' & = & Z
 \end{array}.
\right.
\]
Then the point $(x_0, y_0, 1)$ goes to $(0, 0, 1)$
in these new coordinates.
Let $x = X'/Z'$ and  $y = Y'/Z'$
be affine coordinates of $\mathbb{P}^2 \setminus \{Z'=0 \}$
with respect to these new coordinates,
and $O=(0,0,1)$.
The equation of $C$ in $x, y$ is $f(x, y)=0$,
where
\begin{align}
 f(x,y) &= (x^q-x)(a_0x+a_1y) + (y^q-y)(b_0x+b_1y) \notag \\
        &= x^q (a_0x+a_1y) +y^q(b_0x+b_1y)
          - (a_0 x^2 +(a_1+b_0)xy + b_1y^2). \label{quadricterms}
\end{align}
Note that the condition of irreducibility of (\ref{additionalcondition}) in 
the original coordinates corresponds to
that of the quadric terms in (\ref{quadricterms}) in these affine coordinates.
Hence $O$ is an ordinary double point.
In fact, the embedded tangent cone at $O$ is defined by 
the quadric terms in (\ref{quadricterms})
which is irreducible over $\mathbb{F}_q$, and then 
it splits into two $\mathbb{F}_q$-conjugate lines over $\mathbb{F}_{q^2}$.

Let us investigate the intersection of $C$ and any line $l$ passing through $O$
over the algebraic closure of $\mathbb{F}_q$.
If $l= \{x=0 \}$,
then the $y$-coordinates of points of $C \cap l$
is given by $b_1(y^{q+1}-y^2)$.
Since $b_1 \neq 0$
(otherwise the quadric terms in (\ref{quadricterms}) is reducible over
$\mathbb{F}_q$), each point of $C \cap l \setminus \{O\}$ is nonsingular.
Next we choose a line $l_c$ defined by $y=cx$.
Substitute $cx$ for $y$ in $f(x,y)$.
Then
\begin{equation}\label{eqinx}
( (a_0+ca_1) +c^q(b_0+c b_1) )x^{q+1}
  - ( (a_0+ca_1) +c(b_0+c b_1) )x^2 =0.
\end{equation}
If both coefficients of $x^{q+1}$ and $x^2$ are nonzero,
then each point of $C \cap l_c \setminus \{O\}$ is nonsingular.

Suppose that
the coefficient of $x^2$ is $0$, namely
$(a_0+ca_1) +c(b_0+c b_1)=0.$
Then $(1,c)$ is a root of
$a_0 x^2 +(a_1+b_0)xy + b_1y^2 =0$, and vice versa.
Note that this quadric is the lowest terms of
the local equation
(\ref{quadricterms}) of $C$ around $O$.
Hence $(a_0+ca_1) +c(b_0+c b_1)=0$ if and only if
$l_c$ is the embedded tangent line of one of the two branches
of the singular points $O$.
In this case, $c \not\in \mathbb{F}_q$
because the quadric is irreducible over $\mathbb{F}_q$.
So $(a_0+ca_1) +c^q(b_0+c b_1) \neq 0$.
(Actually, $(a_0+ca_1) +c^q(b_0+c b_1) = 0$ and $(a_0+ca_1) +c(b_0+c b_1)=0$
imply $(c^q-c)(b_0 + cb_1)=0$. Since $b_1$ and $b_0$ 
are in $\mathbb{F}_q$
and at least one of them is nonzero,
$c \in \mathbb{F}_q$, which is a contradiction.)
Hence $C \cap l_c = \{ O \}$.

Finally suppose that
\begin{equation}\label{lastcase}
(a_0+ca_1) +c^q(b_0+c b_1)=0.
\end{equation}
Then $l_c$ meets with $C\setminus \{ O \}$ at $(1, c, 0)$
of multiplicity $q-1$.
Choose local coordinates around $(1, c, 0)$
as $t := \frac{Y'}{X'}-c$ and $z:=\frac{Z'}{X'}$.
Then the local equation of $C$ is 
\begin{align*}
0 &= (1-z^{q-1})(a_0 + a_1(c+t))
    + (c^q + t^q -(c+t)z^q) (b_0 + b_1(c+t)) \\
  &= (a_1 + b_1c^q)t + \text{ (higher terms in $t$ and $z$).}
\end{align*}
Here, if $a_1 + b_1c^q =0$,
then $a_0 + b_0c^q =0$
because the assumption (\ref{lastcase}).
This is a contradiction because
$
 \det \left(
        \begin{array}{cc}
         a_0 & a_1 \\
         b_0 & b_1 
        \end{array}
      \right) \neq 0.
$
Hence the coefficient of $t$ in the local equation is nonzero,
and hence $C$ is nonsingular at $(1, c, 0)$.
\qed


\end{document}